\documentclass[12pt]{article}

\usepackage{amsmath, amsfonts, amssymb}
\usepackage{latexsym}
\usepackage{graphicx}
\usepackage{color} 
\usepackage{url}
\usepackage{lscape}
\usepackage{verbatim}

\newtheorem{thm}{Theorem}[section]
 \newtheorem{cor}[thm]{Corollary}
 
 \newtheorem{prop}[thm]{Proposition}
 
 \newtheorem{rem}[thm]{Remark}

\newcommand{\iy}{\infty}

\newcommand{\bZ}{{\mathbf Z}}

\newcommand{\bR}{{\mathbf R}}

\newcommand{\cF}{{\mathcal F}}

\newcommand{\cL}{{\mathcal L}}

\newcommand{\cM}{{\mathcal M}}

\newcommand{\cV}{{\mathcal V}}

\newcommand{\al}{\alpha}
\newcommand{\be}{\beta}
\newcommand{\ga}{\gamma}

\newcommand{\eps}{\varepsilon}

\newcommand{\La}{\Lambda}

\newcommand{\om}{\omega}
\newcommand{\ph}{\varphi}

\newcommand{\si}{\sigma}

\newcommand{\sm}{\setminus}
\newcommand{\n}{\|}

\newcommand{\ti}{\widetilde}

\newcommand{\1}{{\bf 1}}
\newcommand{\2}{{\bf 2}}

\newcommand{\spn}{{\rm span}}
\newcommand{\rk}{{\rm rk}}

\def\GL{\operatorname{GL}}

\begin{document}
\begin{center}
{\Large \bf  Lattices from tight equiangular frames}

\bigskip
{\large Albrecht B\"ottcher, Lenny Fukshansky, Stephan Ramon Garcia,\\[0.5ex]

Hiren Maharaj, Deanna Needell}
\end{center}

\begin{quote}
{\bf Abstract.} We consider the set of all linear combinations with integer coefficients of the
vectors of a unit tight equiangular $(k,n)$ frame and are interested in the question
whether this set is a lattice, that is, a discrete additive subgroup of the $k$-dimensional
Euclidean space. We show that this is not the case if the cosine of the angle of the frame is irrational.
We also prove that the set is a lattice for $n = k+1$ and that there are infinitely many $k$ such
that a lattice emerges for $n = 2k$. We dispose of all cases in dimensions $k$ at most 9.  In particular,
we show that a (7,28) frame generates a strongly eutactic lattice and give
an alternative proof of Roland Bacher's recent observation that this lattice is perfect.


\end{quote}

\section{Introduction}

Let $2 \le k < n$ and
let $G$ be a real $k \times n$ matrix. Denote the columns of $G$ by $f_1, \ldots, f_n$. These
columns or $G$ itself are called a unit tight equiangular $(k,n)$ frame if
$GG'=\ga I$ with $\ga=n/k$ (tightness) and
$G' G=I+(1/\al)C$ with $\al=\sqrt{k(n-1)/(n-k)}$ and a matrix $C$
whose diagonal entries are zero and the other entries of which are $\pm 1$
(property of being equiangular unit vectors). Define
$\La(G)={\rm span}_\bZ\{f_1, \ldots, f_n\}$. Our investigation is motivated by the following question.

\medskip
{\em When is $\La(G)$ a lattice, that is, a discrete additive subgroup of $\bR^k$? In case it is a lattice, what are its geometric properties?}

\medskip
After having posed the question in its most concise form, some comments are in order. By $\bR^k$ we understand the column-wise written Euclidean $\bR^k$ with the usual
scalar product $(\cdot,\cdot)$. The condition $G' G=I+(1/\al)C$ with $C$ as above means that $\n f_j\n =1$ for all $j$
and that $|(f_i,f_j)|=1/\al$ for $i \neq j$. In other words, the vectors $f_j$ are all unit vectors and each pair
of them makes the angle $\ph$ or $\pi-\phi$ such that $|\cos \ph| = |\cos (\pi-\phi)| =1/\al$.
The equality $GG'=\ga I$ is equivalent to the
requirement that $\n G'x\n^2=\ga \n x\n^2$ for all $x$ in $\bR^k$, which in turn is the same as saying that
$\sum_{j=1}^n (f_j,x)^2=\ga \n x\n^2$ for all $x \in \bR^k$.
It is well known since~\cite{Stroh, Sus} that the two equalities $G' G=I+(1/\al)C$ with $C$ as above and $GG'=\ga I$
necessarily imply that $\ga =n/k$ and
$\al=\sqrt{k(n-1)/(n-k)}$.

\medskip
Tight equiangular frames (TEFs) possess many properties similar to orthonormal bases, yet may also be highly overcomplete, making them very attractive in many applications.  For this reason there has been a recent surge of work addressing the construction and analysis of these frames.  They appear in many practical applications, such as error correcting codes \cite{HP,Stroh}, wireless communications \cite{strohmer2003optimal,Stroh}, security \cite{mixon2011equiangular}, and sparse approximation \cite{fickus2012steiner,rusu2015optimized,tropp2004greed,tsiligianni2012use}.

\medskip
In sparse approximation for example, the incoherence (small $1/\alpha$, the absolute value of pairwise inner products of vectors) of TEFs allows them to be used as sensing operators.  Viewing a TEF as a matrix whose columns consist of the frame vectors, samples of a signal are acquired via inner products between the signal and the rows of this (typically highly underdetermined) matrix.  Under the assumption that the signal vector is sparse (has a small number of nonzero coordinates), the signal can be accurately reconstructed from this compressed representation.  However, in many applications there is more known about the signal than it simply being sparse.  For example, in error correcting codes \cite{candes2005error} and communications applications like MIMO \cite{rossi2014spatial} and cognitive radio \cite{axell2012spectrum}, the signal vectors may come from some lattice.  However, there has been very little rigorous mathematical developments on the intersection between arbitrary lattice-valued signals and sparse approximation (see e.g. \cite{flinth2016promp} and references therein).

\medskip
In this work, we attempt to take the first step toward a rigorous analysis of properties of tight equiangular frames and associated lattices.  We are especially interested in the following questions.  When does the integer span of a TEF form a lattice?  Does this lattice have a basis of minimal vectors?  Is the generating frame contained among the minimal vectors of this lattice?  We also study further geometric properties of the resulting lattices, such as eutaxy and perfection.  Our hope is that this investigation will contribute not only to the understanding of TEFs in general, but also to their explicit use in applications with lattice-valued signals.  For example, if the integer span of a TEF is a lattice, then the image of that TEF viewed as a sensing matrix restricted to integer-valued signals forms a discrete set.  In some sense this is analogous to the well-known Johnson-Lindenstrauss lemma \cite{baraniuk2008simple}, and may be used to provide reconstruction guarantees for TEF sampled signals.  More concretely, if the lattice constructed from the TEF is such that its minimal vectors are the frame vectors themselves, this guarantees a minimum separation between sample vectors in its image.  These types of properties are essential for sparse reconstruction and can be leveraged to design new sampling mechanisms and reconstruction guarantees.  On the other hand, it is also useful to know when such properties are impossible.  We leave a detailed analysis and link to applications as future work, and focus here on the mathematical underpinnings to the questions raised above.

\section{Main results}

Let $\cL$ be a lattice in $\bR^k$, and let $V = \spn_{\bR} \cL$ be the subspace of $\bR^k$ that it spans. Then the rank of $\cL$, denoted by $\rk(\cL)$, is defined to be the dimension of $V$. We say that $\cL$ has full rank if $V=\bR^k$. The minimal distance of a lattice $\cL \subset \bR^k$ is defined as $d(\cL)=\min \{\n x\n: x \in \cL \sm \{0\}\}$. The set of minimal vectors, $S(\cL)$, is the set of all $x \in \cL$ with $\n x\n =d(\cL)$. The lattice $\cL$ is called well-rounded if $\bR^k = \spn_{\bR} S(\cL)$, and we say that it is generated by its minimal vectors if $\cL = \spn_{\bZ} S(\cL)$. It is known that the second condition is strictly stronger than the first when $\rk(\cL) \geq 5$. An even stronger condition (at least when $\rk(\cL) \geq 10$) is that $S(\cL)$ contains a basis for $\cL$, i.e., there exist $\bR$-linearly independent vectors $f_1, \ldots, f_{\rk(\cL)} \in S(\cL)$ such that $\spn_\bZ \{f_1, \ldots, f_{\rk(\cL)}\}=\cL$; if this is the case, we say that $\cL$ has a basis of minimal vectors.

\medskip
A finite subset $\{ q_1,\dots, q_m \}$ of the unit sphere $\Sigma_{k-1}$ in $\bR^k$ is called a spherical $t$-design for a positive integer $t$ if for every real polynomial $p$ of degree $\leq t$ in $k$ variables,
$$\int_{\Sigma_{k-1}} p(x)\ d\si(x) = \frac{1}{m} \sum_{i=1}^k p(q_i),$$
where $d\si$ denotes the unit normalized surface measure on the sphere $\Sigma_{k-1}$. A full rank lattice in $\bR^k$ is called strongly eutactic if its set of minimal vectors (normalized to lie on $\Sigma_{k-1}$) forms a spherical 2-design. We finally define the notion of a perfect lattice. Recall that we write vectors $x$ in $\bR^k$ as column vectors. A full rank lattice $\cL$ in $\bR^k$ is called perfect if the set of  symmetric $k \times k$ matrices $\{ xx' : x \in S(\cL) \}$ spans all  real symmetric $k \times k$ matrices as an $\bR$-vector space.

\medskip
Two lattices $\cL$ and $\cM$ in $\bR^k$ are called similar if $\cL = aU\cM$ for some $a \in \bR$ and some orthogonal $k \times k$ matrix $U$. Conditions such as well-roundedness, generation by minimal vectors, existence of bases of minimal vectors, strong eutaxy, and perfection are preserved on similarity classes of lattices. Furthermore, there are only finitely many strongly eutactic and only finitely many perfect similarity classes of lattices in $\bR^k$ for each $k \geq 1$.

\medskip
Given a full rank lattice $\cL \subset \bR^k$, it is possible to associate a sphere packing to it by taking spheres of radius $d(\cL)/2$ centered at every point of $\cL$. It is clear that no two such spheres will intersect in their interiors. Such sphere packings are usually called lattice packings. One convenient way of thinking of a lattice packing is as follows. The Voronoi cell of $\cL$ is defined to be
$$\cV(\cL) := \{ x \in \bR^k : \| x \| \leq \| x - y\| \ \forall\ y \in \cL \}.$$
Then $\bR^k$ is tiled with translates of the Voronoi cell by points of the lattice, and spheres in the packing associated to $\cL$ are precisely the spheres inscribed in these translated Voronoi cells. A compact measurable subset of $\bR^k$ is called a fundamental domain for a lattice $\cL$ if it is a complete set of coset representatives in the quotient group $\bR^k/\cL$. All fundamental domains of the same lattice have the same volume, and the Voronoi cell of a lattice is an important example of a fundamental domain.

\medskip
A central problem of lattice theory is to find a lattice in each dimension $k \geq 1$ that maximizes the density of the associated lattice packing. There is an easy formula for the packing density of a lattice. A lattice $\cL$ in $\bR^k$ can be written as $\cL = B\bZ^k$, where $B$ is a basis matrix of $\cL$, i.e., the columns of $B$ form a basis for $\cL$. The determinant of $\cL$ is then defined to be $\det \cL := \sqrt{\det (B'B)}$, which is an invariant of the lattice, since any two basis matrices of $\cL$ are related by a integer linear transformation with determinant $\pm 1$. The significance of the determinant is given by the fact that it is equal to the volumes of the fundamental domains. It is then easy to observe that the density of the lattice packing associated to $\cL$ is the volume of one sphere divided by the volume of the translated Voronoi cell that it is inscribed into, that is,
\begin{equation}
\delta(\cL) := \frac{\omega_k d(\cL)^k}{2^k \det \cL}, \label{den}
\end{equation}
where $\omega_k$ is the volume of the unit ball in $\bR^k$. In fact, this packing density function~$\delta$ is defined on similarity classes of lattices in a given dimension, and a great deal of attention in lattice theory is devoted to studying its properties. There is a natural quotient metric topology on the space of all full rank lattices in $\bR^k$, given by identifying this space with $\GL_k(\bR)/\GL_k(\bZ)$: indeed, every $A \in \GL_k(\bR)$ is a basis matrix of some lattice, and $A,B \in \GL_k(\bR)$ are basis matrices for the same lattice if and only if $A=UB$ for some $U \in \GL_k(\bZ)$. A lattice is called extreme if it is a local maximum of the packing density function in its dimension: this is a particularly important class of lattices that are actively studied. A classical result of Voronoi states that perfect strongly eutactic lattices are extreme (see, for instance, Theorem 4 of~\cite{achill1}); on the other hand, if a lattice is strongly eutactic, but not perfect, then it is a local minimum of the packing density function (see Theorem~9.4.1 of~\cite{martinet}). A good source for further information about lattice theory is Martinet's book~\cite{martinet}.

\medskip
We now return to our construction $\La(G)$ from unit equiangular frames and describe our results.
It is well known that unit tight equiangular $(k,k+1)$ frames exist for all $k \ge 2$.
According to \cite{Sus}, except for the $(k,k+1)$-case, the only unit tight equiangular $(k,n)$ frames
with $k \le 9$ are
\begin{equation}
(3,6), (5,10), (6,16), (7,14), (7,28), (9,18)\label{flist}
\end{equation}
frames. Our first result says the following.

\medskip
\begin{prop} \label{Prop I1}
If $\La(G)$ is a lattice, then $\al$ must be a rational number.
\end{prop}

Thus, since $\al=1/\sqrt{5}$ for the $(3,6)$ frame, $\al =1/\sqrt{13}$ for the $(7,14)$ frame,
and $\al=1/\sqrt{17}$ for the $(9,18)$ frame,
these three frames do not generate lattices. We will show that there are unit tight equiangular
$(5,10)$, $(6,16)$, and $(7,28)$ frames which generate lattices.
Moreover, we will prove the following results.

\begin{thm} \label{Theo I2}
{\rm (a)} For every $k \ge 2$, there are unit tight equiangular $(k,k+1)$ frames $G$ such that $\La(G)$
is a full rank lattice. The lattice $\La(G)$  has a basis of minimal vectors, it is non-perfect and strongly eutactic, and hence it is a local minimum of the packing density function in dimension $k$.
\

\medskip
{\rm (b)} There are infinitely many $k$ for which there exist unit tight equiangular $(k,2k)$ frames
$G$ such that $\La(G)$ is a full rank lattice.

\medskip
{\rm (c)} There is a unit tight equiangular $(7,28)$ frame $G$ for which $\La(G)$ has a basis of minimal vectors, is a perfect strongly eutactic lattice, and hence extreme.
\end{thm}

\begin{rem} {\rm We explicitly construct the lattices of Theorem~\ref{Theo I2}. We show that those of parts (a) and (c) and those with $k \leq 13$ of part (b) have the property that the set of minimal vectors consists precisely of $\pm$ the generating frame vectors. The well known result of Gerzon (see, for instance, Theorem~C of~\cite{Sus}) asserts that for a $(k,n)$ tight equiangular frame
necessarily $n \leq k(k+1)/2$. On the other hand, $k(k+1)/2$ is the minimal number of ($\pm$ pairs of) minimal vectors necessary (but not sufficient) for a lattice in~$\bR^k$ to be perfect. Since only very few tight equiangular frames achieve equality in Gerzon's bound, it is likely quite rare for perfect lattices to be generated by tight equiangular frames. Perfection is a necessary condition for extremality, and hence it is unreasonable to expect to obtain extreme lattices often in this way. The only such example we have discovered is the lattice from the $(7,28)$ frame in part (c) of our Theorem~\ref{Theo I2}, perfection of which has also previously been discussed in~\cite{bacher}.}
\end{rem}

The strong eutaxy of our lattice constructions in Theorem~\ref{Theo I2}(a),(c) is established directly with the use of the following result.

\begin{prop} \label{eutaxy} Suppose that $\Lambda(G)$ is a lattice and $S(\Lambda(G)) = \{ \pm f_1,\dots,\pm f_n \}$. Then $\Lambda(G)$ is strongly eutactic.
\end{prop}

{\em Proof.}
A spanning set $\{ g_1,\dots,g_m \}$ for $\bR^k$ is called a Parseval frame if
$\|x\|^2 = \sum_{j=1}^m (g_j,x)^2$
for all $x \in \bR^k$. Further, $\{ g_1,\dots,g_m \}$ is a spherical 2-design if and only if
$$\left\{ \sqrt{k/m}\ g_1,\dots,\sqrt{k/m}\ g_m \right\}$$
is a Parseval frame and $\sum_{i=1}^m g_i = 0$ (see~\cite{HP} for details, especially Proposition 1.2).

\medskip
Now let $G = (f_1\ \dots\ f_n )$ be a unit tight equiangular $(k,n)$ frame, and assume that $\Lambda(G)$ is a lattice such that $S(\Lambda(G)) = \{ \pm f_1,\dots,\pm f_n \}$. We
then have
\begin{eqnarray*}
\| x \|^2 & = & \frac{k}{2n} \sum_{j=1}^n \left( ( x, f_j )^2 + ( x, -f_j )^2 \right) \\
& = & \sum_{j=1}^n \left( \left( x, \sqrt{\frac{k}{2n}} f_j \right)^2 + \left( x, -\sqrt{\frac{k}{2n}} f_j \right)^2 \right),
\end{eqnarray*}
for every $x \in \bR^k$. Hence $\left\{ \pm \sqrt{k/2n}\ f_1,\dots,\pm \sqrt{k/2n}\ f_n \right\}$ is a Parseval frame, and therefore $S(\Lambda(G))$ is a spherical 2-design.
$\;\:\square$

\medskip
A summary of a part of our results is given in Table~\ref{table1}.

\begin{table}
\caption{Summary of a part of our results.}
\label{table1}
\[\begin{array}{|l|l|l|l|}
\hline
~ & ~ & ~ & ~\\
(k,n) & \mbox{cosine} \;{\displaystyle \frac{1}{\al}} & \mbox{Volume of a} & S(\La)=\{\pm f_1, \ldots, \pm f_n\}?\\
& & \mbox{fundamental domain}  & \mbox{Basis of minimal vectors?}\\
~ & ~ & ~ & ~\\
\hline
\hline
~ & ~ & ~ & ~\\
(k+1,k) & {\displaystyle \frac{1}{k}} & {\displaystyle\frac{1}{\sqrt{k+1}}\left(1+\frac{1}{k}\right)^{k/2}}
& \mbox{Yes, Yes} \\
~ & ~& ~ & ~\\
\hline
~ & ~ & ~ & ~\\
(3,6) & {\displaystyle \frac{1}{\sqrt{5}}=0.4472} & \mbox{no lattice} & \\
~ & ~ & ~ & ~\\
\hline
~ & ~ & ~ & ~\\
(5,10) & {\displaystyle\frac{1}{3}} & {\displaystyle  \frac{4}{9}=0.4444}
& \mbox{Yes, Yes} \\
~ & ~ & ~ & ~\\
\hline
~ & ~ & ~ & ~\\
(6,16) & {\displaystyle \frac{1}{3}} & {\displaystyle\frac{2^3}{3^3}=0.2963}
& \mbox{Yes, Yes} \\
~ & ~ & ~ & ~\\
\hline
~ & ~ & ~ & ~\\
(7,14) & {\displaystyle \frac{1}{\sqrt{13}}=0.2774} & \mbox{no lattice} & \\
~ & ~ & ~ & ~\\
\hline
~ & ~ & ~ & ~\\
(7,28) & {\displaystyle \frac{1}{3}} & {\displaystyle \frac{2^3}{3^{7/2}}=0.1711}
& \mbox{Yes, Yes, and perfect} \\
~ & ~ & ~ & ~\\
\hline
~ & ~ & ~ & ~\\
(9,18) & {\displaystyle \frac{1}{\sqrt{17}}=0.2425} & \mbox{no lattice} & \\
~ & ~ & ~ & ~\\
\hline
~ & ~ & ~ & ~\\
(13,26) & {\displaystyle \frac{1}{5}} & {\displaystyle \frac{2^6}{5^{9/2}}=0.0458} & \mbox{Yes, Yes}\\
~ & ~ & ~ & ~\\
\hline
~ & ~ & ~ & ~\\
(25,50) & {\displaystyle \frac{1}{7}} & {\displaystyle \frac{2^{11}\cdot 3 \cdot 5 \cdot 11^2}{7^{23/2}}=0.00071052} & \mbox{?, ?}\\
~ & ~ & ~ & ~\\
\hline
\end{array}
\]
\end{table}

\section{Rationality of the cosine of the frame}

Suppose $G$ is a unit tight $(k,n)$ frame. Then $GG'=\ga I$ and hence $G$ has rank $k$.
Let $G_0$ be the $k \times k$ matrix formed by arbitrarily chosen $k$ linearly independent
columns of $G$ and denote by $G_1$ the $k \times (n-k)$ matrix constituted by the
remaining columns. We may without loss of generality assume that $G=(G_0\; G_1)$.
We emphasize that $G_0$ is invertible. Recall that $\La(G)$ is called a full-rank lattice if
${\rm span}_\bR\{f_1, \ldots, f_n\}$ is all of $\bR^k$. Note that in the following proposition
we do not require equiangularity.

\begin{prop} \label{Prop 1}
Let $G=(G_0 \; G_1)$ be a unit tight $(k,n)$ frame. Then the following are
equivalent.

\medskip
(i) $\La(G)$ is a lattice.

\medskip
(ii) $\La(G)$ is a full rank lattice.

\medskip
(iii) There exist $\be \in \bZ\sm\{0\}$ and $X\in \bZ^{k \times (n-k)}$ such that
$G_0^{-1}G_1=(1/\be)X$.

\medskip
If (iii) holds with $\be=1$, then $G_0$ is a basis matrix for $\La(G)$.
\end{prop}

{\em Proof.} Since $G_0$ is invertible, we have ${\rm span}_\bR\{f_1,\ldots,f_n\}=\bR^k$,
which proves the equivalence of (i) and (ii). Suppose (ii) holds. Then $G_0=BX_0$
and $G_1=BX_1$ with
an invertible $k \times k$ matrix $B$ and integer matrices $X_0,X_1$. The matrix $X_0$ is invertible,
so $B=G_0X_0^{-1}$ and hence
\[G_1=G_0X_0^{-1}X_1=G_0\frac{1}{\det X_0}X_2X_1=\frac{1}{\be}G_0X\]
with $\be =\det X_0$ and $X=X_2X_1$. This proves (iii). Conversely, suppose (iii)
is true. It is clear that $\La(G)={\rm span}_\bZ\{f_1,\ldots,f_n\}$ is an additive subgroup
of $\bR^k$. Put $B=(1/\be)G_0$. Then $B$ is invertible, $G_0=B X_0$ with $X_0=\be I$
and $G_1=BX_1$ with $X_1=X$. It follows that $\La(G)$
is a subset of $L_B:=\{BZ: Z \in \bZ^{k \times 1}\}$. As the latter set is discrete, so must be
$\La(G)$. This  proves (i). Finally, if $\be=1$, then $B=G_0$, which implies that
$L_B \subset \La(G)$ and hence $L_B=\La(G)$. Consequently, $B$ is a basis matrix for $\La(G)$. $\;\:\square$

\begin{prop} \label{Prop 2}
Let $G=(G_0 \; G_1)$ be a unit tight equiangular $(k,n)$ frame. If $\La(G)$ is a lattice,
then $\al$ must be a rational number.
\end{prop}

{\em Proof.} By Proposition \ref{Prop 1}, we may assume that $\La(G)$ is a full rank
lattice. So $G=BZ$ with an invertible matrix $B$ and a matrix $Z \in \bZ^{k \times n}$.
Multiplying the equality $\ga I =GG'=BZZ'B'$ from the right by $(B')^{-1}$ and then from
the left by $B'$, we obtain $\ga I=B'B ZZ'$ and thus,
\[(I+(1/\al)C)Z'=G'GZ'=Z'B'BZZ'=\ga Z',\]
which implies that $CZ'=\al(\ga-1)Z'$. If $\al$ is irrational, the last equality yields $Z=0$,
and this gives $G=0$, a contradiction.  $\;\:\square$

\bigskip
The previous proposition implies in particular that a unit tight equiangular $(3,6)$ frame does not induce a lattice.
The reader might enjoy to see the reason for this failure also from the following perspective.
Consider the tight unit equiangular $(3,6)$ frame $G$ that is
induced by the $6$ upper vertices of a regular icosahedron. As shown in~\cite{Sus}, with $p=(1+\sqrt{5})/2$, this frame is given by the columns of the matrix
\[G=\frac{1}{\sqrt{1+p^2}}
\left(\begin{array}{rrrrrr}
0 & 0 & 1 & -1 & p & p\\
1 & -1 & p & p & 0 & 0 \\
p & p & 0 & 0 & 1 & -1\end{array}\right).\]
We have $c =1/\sqrt{5}$.
By Dirichlet's approximation theorem, there are integers $x_n,y_n$ such that $y_n \to \iy$ and
\[\left|\frac{x_n}{y_n}+p\right| \le \frac{1}{y_n^2}.\]
In particular, $x_n+py_n \to 0$ as $n \to \iy$
The linear combination of the columns of $G$ with the coefficients $x_n+y_n, \:y_n-x_n,\: y_n,\: y_n, \: x_n, \: -x_n$,
equals
\[\frac{1}{\sqrt{1+p^2}}\left(\begin{array}{c} 0 \\ 2(x_n+y_n p)\\2(y_n p +x_n) \end{array}\right),\]
which tends to zero as $n \to \iy$. Consequently, $\La(G)$ is not a discrete subgroup of $\bR^3$ and thus
it is not a lattice.

\section{Unit tight equiangular $\boldsymbol{(k,2k)}$ frames}

We first consider the case $n=2k$. Then $\ga =2$ and $\al=\sqrt{n-1}$. We furthermore suppose
that $n=p^r+1$ with an odd prime number $p$ and a natural number $r$. If $r$ is odd and
$p=4\ell+3$, then $k$ is even, which implies that unit tight equiangular $(k,n)$ frames
do not exist (Theorem~17 of~\cite{Sus}). If $r$ is odd and $p=4\ell+1$, then
unit tight equiangular $(k,n)$ frames $G$ exist, but $\La(G)$ is not a lattice because
$\al$ is irrational. We are so left with the case where $r$ is even.

\begin{thm} \label{Theo 1}
Let $ k \ge 2$ and $n=2k$. If $n=p^{2m}+1$ with an odd prime number $p$ and a natural number $m$,
then there exists a unit tight equiangular $(k,n)$ frame $G$ such that $\La(G)$ is a full
rank lattice.
\end{thm}

{\em Comments.} This theorem proves Theorem~\ref{Theo I2}(b) and will be a consequence of the following Theorem~\ref{Theo 1a}.
Before turning to the proof of Theorem~\ref{Theo 1a}, which is a combination of ideas of Goethals
and Seidel~\cite{Goe} and Strohmer and Heath~\cite{Stroh}, some comments seem to be in order. Following~\cite{Stroh},
we start with a symmetric $n \times n$ conference matrix $C$, that is, with a symmetric matrix $C$ that has zeros on the
main diagonal and $\pm 1$ elsewhere and that satisfies $C^2=(n-1)I$. Under the hypothesis of Theorem~\ref{Theo 1},
such matrices were first constructed by Paley~\cite{Pa}. Goethals and Seidel~\cite{Goe} showed that one can always
obtain such matrices in the form
\begin{equation}
C=\left(\begin{array}{rr} A & D \\ D & -A \end{array}\right)\label{CZirk}
\end{equation}
where $A$ and $D$ are symmetric $k \times k$ circulant matrices. Let $a$ and $b$ be any rational numbers such that
$a^2+b^2=\al^2$ ($=n-1=p^{2m}$). Theorem 3.4 of~\cite{Goe} says that, under certain conditions, one
can in turn represent the matrix~(\ref{CZirk})
as
\begin{equation}
\left(\begin{array}{rr} A & D \\ D & -A \end{array}\right)
=\left(\begin{array}{rr} I & -N \\ N & I \end{array}\right)^{-1}
\left(\begin{array}{rr} aI & bI  \\ bI  & -aI\end{array}\right)
\left(\begin{array}{rr} I & -N \\ N & I \end{array}\right)\label{Nab}
\end{equation}
with a symmetric circulant matrix $N$ all entries of which are rational numbers. The conditions ensuring the
representation~(\ref{Nab}) are that $D+bI$ or $A+aI$ are invertible. We have
\begin{equation}
 N=(A+aI)^{-1}(bI-D) \quad \mbox{or} \quad N=(D+bI)^{-1}(A-aI)\label{NDef}
\end{equation}
if $A+aI$ or $D+bI$ is invertible, respectively. (Note that all occurring blocks are symmetric
circulant matrices and in particular commuting matrices.)
As there are infinitely many different
decompositions of $p^{2m}$ into the sum of two squares of rationals,
for example,
\[p^{2m}=\left(\frac{t^2-s^2}{t^2+s^2}p^m\right)^2+\left(\frac{2ts}{t^2+s^2}p^m\right)^2\]
with integers $s$ and $t$, we can, for given $A$ and $D$, always find rational $a$ and $b$ such that $a^2+b^2=\al^2$
and both $D+bI$ and $A+aI$ are invertible.

\medskip
Let, for example $n=10$. A matrix (\ref{CZirk}) with symmetric circulant matrices $A$ and $D$ is completely given by
its first line, which is of the form
\[0,\eps_1,\eps_2,\eps_2,\eps_1,\quad \eps_3, \eps_4,\eps_5,\eps_5,\eps_4\]
with $\eps_j \in \{-1,1\}=:\{-,+\}$. These are $2^5=32$ matrices. Exactly four of them satisfy $C^2=9I$. Their first lines
and the eigenvalues of $D$ are
\begin{eqnarray}
& &  0,    -,     +,     +,    -, \quad  -,     +,     +,     +,     +, \qquad -2,-2,-2,-2,3\label{t1} \\
& &  0,    -,     +,     +,    -,  \quad    +,    -,    -,    -,    -, \qquad -3, 2,2,2,2,\label{t2}\\
& &  0,     +,    -,    -,     +,  \quad  -,     +,     +,     +,     +, \qquad -2,-2,-2,-2,3 \label{t3}\\
& &  0,     +,    -,    -,     +,   \quad  +,    -,    -,    -,    -, \qquad -3,2,2,2,2.\label{t4}
\end{eqnarray}
The corresponding matrix $A$ is always singular. We see that in all cases we may
take $a=3$ and $b=0$ ($3^2+0^2=9$) because $D$ is invertible. In the cases~(\ref{t1}) and~(\ref{t3})
we could also take $a=0$ and $b=3$ ($0^2+3^2=9$) since $D+3I$ is invertible. In fact, we will prove
the following theorem. As shown above, the hypothesis of this theorem can always be satisfied,
so that this theorem implies Theorem~\ref{Theo 1}.

\begin{thm} \label{Theo 1a}
Let $ k \ge 2$ and $n=2k$. Suppose $n=p^{2m}+1$ with an odd prime number $p$ and a natural number $m$,
let $a$ and $b$ be rational numbers such that $a^2+b^2=p^{2m}$ and $a \neq -p^m$. Let $A$ and $D$ be symmetric
$k \times k$ circulant
matrices such that~(\ref{CZirk}) is a conference matrix, and assume $A+aI$
or $D+bI$ is invertible. Define $N$ by~(\ref{NDef}) and put $\al=\sqrt{n-1}=p^m$. Then
\begin{equation}
G=\frac{1}{\sqrt{\al(\al+a)}}(I+N^2)^{-1/2}\left[\begin{array}{cc}
(\al+a)I+bN & bI-(\al+a)N\end{array}\right]\label{Gnew}
\end{equation}
is a unit tight equiangular $(k,n)$ frame $G$ such that $\La(G)$ is a full
rank lattice.
\end{thm}

{\em Proof of Theorem \ref{Theo 1a}.} The requirement $a \neq -p^m$ assures that $\al+a\neq 0$. Let
\[W=\left(\begin{array}{ll}W_{11} & W_{12} \\ W_{21} & W_{22}\end{array}\right)
=\frac{1}{\sqrt{2\al(\al+a)}}(I+N^2)^{-1/2}\left(\begin{array}{ll} U_{11} & U_{12} \\ U_{21} & U_{22}\end{array}\right)\]
with
\[\left(\begin{array}{ll} U_{11} & U_{12} \\ U_{21} & U_{22}\end{array}\right)
=\left(\begin{array}{cc}
(\al+a)I+bN & bI-(\al+a)N\\ bI-(\al+a)N & -\al I- (\al+a)N\end{array}\right).\]
Using (\ref{Nab}) one can show by straightforward computation that
\[C\left(\begin{array}{ll} U_{11} \\ U_{21}\end{array}\right)
=\al \left(\begin{array}{ll} U_{11} \\ U_{21}\end{array}\right),
\quad
C\left(\begin{array}{ll} U_{12} \\ U_{22}\end{array}\right)
=-\al \left(\begin{array}{ll} U_{12} \\ U_{22}\end{array}\right),\]
which implies that
\[C\left(\begin{array}{ll} U_{11} & U_{12} \\ U_{21} & U_{22}\end{array}\right)
= \left(\begin{array}{ll} U_{11} & U_{12} \\ U_{21} & U_{22}\end{array}\right)
\left(\begin{array}{rr} \al I & 0 \\ 0 & -\al I\end{array}\right)\]
and thus
\begin{equation}
C(I+N^2)^{1/2}W=(I+N^2)^{1/2}W\left(\begin{array}{rr} \al I & 0 \\ 0 & -\al I\end{array}\right).\label{CNW}
\end{equation}
We have $W^2=I$. Indeed,
\[U_{11}^2+U_{12}U_{21}=U_{21}U_{12}+U_{22}^2=[(\al+a)^2+b^2](I+N^2)=2\al(\al+a)(I+N^2),\]
whence $W_{11}^2+W_{12}W_{21}=W_{21}W_{12}+W_{22}^2=I$, and similarly one gets
that the off-diagonal blocks of $W^2$ are zero. From~(\ref{CNW}) we therefore get
\[
C  =  (I+N^2)^{1/2}W\left(\begin{array}{rr} \al I & 0 \\ 0 & -\al I\end{array}\right)W(I+N^2)^{-1/2}
\left(\begin{array}{cc} I & 0 \\ 0 & I\end{array}\right)
 =  W\left(\begin{array}{rr} \al I & 0 \\ 0 & -\al I\end{array}\right)W,
\]
or equivalently,
\begin{equation}
I+\frac{1}{\al}C=W\left(\begin{array}{rr} 2 I & 0 \\ 0 & 0\end{array}\right)W
=2 \left(\begin{array}{cc} W_{11}^2 & W_{11}W_{12} \\ W_{21}W_{11} & W_{21}W_{12}\end{array}\right)\label{WW1}.
\end{equation}
The matrix $G$ given by (\ref{Gnew}) is just $G=\sqrt{2}(W_{11}\;W_{21})$. We claim that $G$ is a unit tight equiangular $(k,n)$ frame.
First, since $W_{11}$ and $W_{21}$ are symmetric, we have
\begin{equation}
G'G = 2\left(\begin{array}{l} W_{11}\\ W_{21}\end{array}\right)
\left(\begin{array}{cc} W_{11} & W_{21}\end{array}\right)
=2\left(\begin{array}{cc} W_{11}^2 & W_{11}W_{21} \\ W_{21}W_{11} & W_{21}^2\end{array}\right), \label{WW2}
\end{equation}
and since $W_{21}=W_{12}$, the right-hand sides of (\ref{WW1}) and (\ref{WW2}) coincide. This proves
that $G$ is unit and equiangular. Secondly,
\[
GG'  =  2 \left(\begin{array}{cc} W_{11} & W_{21}\end{array}\right)\left(\begin{array}{l} W_{11}\\ W_{21}\end{array}\right)
=2 (W_{11}^2 +W_{21}^2)=2I,
\]
which shows that $G$ is tight with $\ga = 2 = n/k$. The equality $GG'=2I$ implies that the rank of $G$ is $k$.
Thus, $G=\sqrt{2}(W_{11} \; W_{12})$ has $k$ linearly independent columns.
We permute the columns of $G$ so that these $k$ linearly independent columns become the first $k$ columns. The
resulting matrix, which is anew denoted by $G$, is a unit tight equiangular $(k,n)$ frame of the form
$G=(G_0 \;G_1)$ with an invertible matrix $G_0$. Furthermore, we have $G_0=(I+N^2)^{-1/2}R$
and $G_1=(I+N^2)^{-1/2}S$ with matrices $R$ and $S$ whose entries are rational numbers. We therefore obtain that
$G_0^{-1}G_1=R^{-1}S$ is a matrix with rational entries, and hence, by Proposition~\ref{Prop 1}, the set $\La(G)$
is a full rank lattice. $\;\:\square$

\begin{cor} \label{Cor 1}
Let $ k \ge 2$ and $n=2k$. Suppose $n=p^{2m}+1$ with an odd prime number $p$ and a natural number $m$,
let $A$ and $D$ be symmetric $k \times k$ circulant matrices such that the matrix (\ref{CZirk}) is
a conference matrix. Put $\al=\sqrt{n-1}=p^m$.
If the matrix $D$ is invertible, then $I\pm(1/\al)A$ are positive definite matrices and, with
the invertible matrix $N:=D^{-1}(A-\al I)$,
\[G:=\sqrt{2}(I+N^2)^{-1/2}\left(\begin{array}{rr} I & -N\end{array}\right)\]
is a unit tight equiangular $(k,n)$ frame $G$ and the set $\La(G)$ is a full
rank lattice. If $N \in \bZ^{k \times k}$, then $G$ may be written as
\begin{equation}
G=B_+\left(\begin{array}{rr} I & -N\end{array}\right)\quad\mbox{with}\quad B_+:=\sqrt{I+(1/\al)A},\label{BP}
\end{equation}
and $B_+$ is a basis matrix for $\La(G)$, while if
$N^{-1} \in \bZ^{k \times k}$, then $G$ may be written in the form
\begin{equation}
G=B_-S\left(\begin{array}{rr} -N^{-1} & I\end{array}\right)\quad\mbox{with}\quad B_-:=\sqrt{I-(1/\al)A},\label{BM}
\end{equation}
where $S:=D|D|^{-1}$ and $|D|$ is the positive definite square root of $D'D$,
and this time $B_-S$ is a basis matrix for $\La(G)$,  Furthermore,
\[\det B_\pm=\sqrt{\det\left(I\pm\frac{1}{\al}A\right)}=\frac{1}{\al^{k/2}}\sqrt{\det (\al I\pm A)}.\]
\end{cor}

\medskip
{\em Remark.} Recall that the determinant ($=$ volume of a fundamental domain) of a lattice is defined as the square root
of $\det (B'B)$ where $B$ is any basis matrix. Thus, if $N$ is an integer matrix, then the determinant of
the lattice is simply $\sqrt{\det(B_+'B_+)}=\det B_+$, while if $N^{-1}$ has integer entries, the determinant
of the lattice $\La(G)$ is
\[\sqrt{\det(S'B_-'B_-S)}=\sqrt{\det (SB_-'B_-S)}=\sqrt{\det (B_-'B_-S^2)}=\det B_-\]
because $S=S'$ and $S^2=I$.

\medskip
{\em Proof.} Since $D$ is invertible, we may use Theorem \ref{Theo 1a} with $a=\al$ and $b=0$ ($\al^2+0^2=\al^2$)
and with $N=D^{-1}(A-\al I)$. In this special case, formula~(\ref{Gnew}) becomes
\begin{equation}
G=\sqrt{2}(I+N^2)^{-1/2}\left(\begin{array}{rr} I & -N\end{array}\right),\label{GNN}
\end{equation}
and since $N$ has rational entries,  Proposition~\ref{Prop 1} implies that $\La(G)$ is a full rank lattice.
Proposition~\ref{Prop 1} also shows that $\sqrt{2}(I+N^2)^{-1/2}$ is a basis matrix for the lattice provided
$N \in \bZ^{k \times k}$. Writing~(\ref{GNN}) as
\[G=-\sqrt{2}(I+N^2)^{-1/2}N\left(\begin{array}{rr} -N^{-1} & I\end{array}\right)\]
and permuting $(-N^{-1}\;I)$ to $(I \;-N^{-1})$ we can deduce from Proposition~\ref{Prop 1} that the matrix $-\sqrt{2}(I+N^2)^{-1/2}N$ is a basis matrix
provided $N^{-1} \in \bZ^{k \times k}$. It remains to show that these two basis matrices are just
the matrices $B_\pm$.

\medskip
As the square of the matrix~(\ref{CZirk}) is $\al^2 I$, we have $A^2+D^2=\al^2 I$.
Since $0$ is not in the spectrum of $D$, the equality $A^2+D^2=\al^2 I$ implies that the spectrum ($=$ set of eigenvalues) of
$A$ is contained in the open interval $(-\al,\al)$. Hence $\al I \pm A$ are positive definite.
Moreover, we
get $D^2=\al^2I -A^2=(\al I -A)(\al I+A)$, and
since all involved matrices are circulants and therefore commute,
we obtain
\begin{eqnarray*}
I+ N^2 & = & I + D^{-2}(A-\al I)^2 = I+D^{-2}(\al I-A)(\al I-A)\\
& = & I +(\al I-A)^{-1}(\al I+A)=  (\al I+A)^{-1}[\al I+A+\al I-A]\\
& = & 2\al(\al I+A)^{-1}=2(I+(1/\al)A)^{-1}.
\end{eqnarray*}
Consequently,
$\sqrt{2}(I+N^2)^{-1/2}=(I+(1/\al)A)^{1/2}=B_+$, which proves~(\ref{BP}).
The matrix $|D|$ is again a circulant matrix and we
have $D=S|D|$ with a circulant matrix $S$ satisfying $S^2=I$. From the equality
$D^2=(\al I-A)(\al I+A)$ we obtain that $|D|=(\al I -A)^{1/2}(\al I+A)^{1/2}$. Thus,
\begin{eqnarray*}
 -\sqrt{2}(I+N^2)^{-1/2}N & = & (I+(1/\al)A)^{1/2} D^{-1}(\al I-A)\\
& = &  \frac{1}{\sqrt{\al}}(\al I+A)^{1/2}(\al I-A)^{-1/2}(\al I+A)^{-1/2}(\al I-A)S\\
& =  &  \frac{1}{\sqrt{\al}}(\al I-A)^{1/2}S=(I-(1/\al)A)^{1/2}S=B_-S.
\end{eqnarray*}
This proves (\ref{BM}). The determinant formulas are obvious.
$\;\:\square$

\bigskip
{\bf Two lattices from  (5,10) frames.}
Let $(A,D)$ be one of the four pairs given by~(\ref{t1}) to~(\ref{t4}). Thus, $k=5$, $n=10$, $\al=3$.
In either case, $D$ is invertible with $\det D=\pm 48$  and we have $\det (3I+ A)=48$. (The eigenvalues of
$A$ are $-\sqrt{5},-\sqrt{5},0,\sqrt{5},\sqrt{5}$.) The four circulant matrices $N=D^{-1}(A-3I)$ have the first rows
\[(+,0,-,-,0),\quad (-,0,+,+,0), \quad (+,-,0,0,-), \quad (-,+,0,0,+).\]
Thus, $N \in \bZ^{5 \times 5}$, and so
by Corollary~\ref{Cor 1}, $\La(G)$ is a lattice, $B_+=\sqrt{I+(1/3)A}$ is a basis matrix, and
$\det B_+=3^{-5/2}\sqrt{48}=2^2/3^2=0.4444\ldots$. (Incidentally, the matrices $N^{-1}$ also have integer entries.) The eigenvalues of $I+(1/3)A$ are
\[1-\frac{1}{3}\sqrt{5}, \quad 1-\frac{1}{3}\sqrt{5}, \quad 1, \quad 1+\frac{1}{3}\sqrt{5}, \quad 1+\frac{1}{3}\sqrt{5}, \quad\]
and hence the smallest eigenvalue of $B_+$ is
$\sqrt{1-(1/3)\sqrt{5}} = 0.50462... >1/2$.
We have $B_+=U E U'$ with an orthogonal matrix $U$ and the diagonal matrix $E$ of the eigenvalues of $B_+$.
Consequently,
\[\n B_+x\n^2=\n UEU'x\n^2=\n EU'x\n^2 >\frac{1}{4}\n U'x\n^2 =\frac{1}{4}\n x\n^2,\]
and hence $\n B_+x\n^2 > 1$
if $\n x\n^2 \ge 4$. So consider the $x \in \bZ^5\sm \{0\}$ with $\n x\n^2 \le 3$. Such $x$ contain only $0,+1,-1$,
and using Matlab we checked that $\n B_+x\n^2 <1.1$ for exactly $20$ nonzero $x$ of these $3^5-1=242$
possible $x$. The $20$ columns $B_+x$ are just $\pm$ the columns of $G$. Thus, {\em the minimal distance of $\La(G)$ is $1$,
$\La(G)$ has a basis of minimal vectors,
and $S(\La(G))=\{\pm f_1, \ldots, \pm f_{10}\}$.} Note that if we denote the basis matrices for the lattices corresponding to~(\ref{t1}) and~(\ref{t4}) by $B_1, \ldots, B_4$,
then actually $B_1=B_2$ and $B_3=B_4$. However, $B_1B_3^{-1}$ is not a scalar multiple of an orthogonal matrix.

\medskip
To ``see'' a concrete matrix $B_+$, note that
in the case where the matrices $A,D$ are specified by~(\ref{t1}),
we obtain that $B_1=B_+=\sqrt{I+(1/3)A}$ equals
\begin{verbatim}
    0.9303   -0.1651    0.2000    0.2000   -0.1651
   -0.1651    0.9303   -0.1651    0.2000    0.2000
    0.2000   -0.1651    0.9303   -0.1651    0.2000
    0.2000    0.2000   -0.1651    0.9303   -0.1651
   -0.1651    0.2000    0.2000   -0.1651    0.9303 .
\end{verbatim}
With the Fourier matrix $F_5=(1/\sqrt{5})(\om^{(j-1)(k-1)})_{j,k=1}^5$, $\om =e^{2\pi i/5}$, this
is $B_1=F_5^*EF_5^{}$ with
\begin{eqnarray*}
E & = & {\rm diag}\left(1, \quad \sqrt{1-\frac{1}{3}\sqrt{5}}, \quad \sqrt{1+\frac{1}{3}\sqrt{5}},
\quad \sqrt{1+\frac{1}{3}\sqrt{5}}, \quad \sqrt{1-\frac{1}{3}\sqrt{5}}\right)\\
& = & {\rm diag}\left(1, \quad \frac{\sqrt{5}-1}{\sqrt{6}}, \quad \frac{\sqrt{5}+1}{\sqrt{6}},
\quad \frac{\sqrt{5}+1}{\sqrt{6}}, \quad \frac{\sqrt{5}-1}{\sqrt{6}} \right).
\quad \square
\end{eqnarray*}
It follows in particular that the numerical values shown above are ${\tt 0.2000}=1/5$, ${\tt 0.9303}=1/5+2 \sqrt{2/15}$,
${\tt -0.1651}=1/5-\sqrt{2/15}$.

\bigskip
{\bf Three lattices from  (13,26) frames.} Let now $k=13$, $n=26$, $\al=5$.
Let $A$ and $D$ be symmetric $13 \times 13$ circulant matrices whose first rows are
\[(0, \eps_1,\eps_2,\eps_3,\eps_4,\eps_5,\eps_6,\eps_6,\eps_5,\eps_4,\eps_3,\eps_2,\eps_1)\]
and
\[(\eps_7, \eps_8,\eps_9,\eps_{10},\eps_{11},\eps_{12},\eps_{13},\eps_{13},\eps_{12},\eps_{11},\eps_{10},\eps_9,\eps_8)\]
with $\eps_k \in\{-1,+1\}=:\{-,+\}$, respectively. There are $2^{13}=8\,192$ such matrices. For exactly $12$ of them the
matrix~(\ref{CZirk}) satisfies $C^2=25I$. The determinant of $D$ is always $\det D=\pm 768\,000= \pm 2^{12}\cdot 3\cdot 5^4$.
Thus, by Corollary~\ref{Cor 1}, $\La(G)$ is a full rank lattice.

\medskip
In exactly $6$ cases, for example if the first rows of $A$ and $D$ are
\[(0,    -,    -,    -,     +,    -,     +,     +,    -,     +,    -,    -,    -) \quad \mbox{and}\quad
(-,    -,     +,     +,     +,    -,     +,     +,    -,     +,     +,     +,    -),\]
we have $N \in \bZ^{13 \times 13}$. We denote the $A$ and $B_+=\sqrt{I+(1/5)A}$ corresponding to these cases
by $A_1, \ldots, A_6$ and $B_1, \ldots, B_6$. Corollary~\ref{Cor 1} implies that $B_j$ is a basis matrix for
the $j$th lattice.
We have
$\det(5I+A_j)=2\,560\,000=2^{12} \cdot 5^4$ and hence $\det B_j=2^6/5^{9/2} \approx 0.0458$ for $1 \le j \le 6$.
In the other $6$ cases, for instance if
the first rows of $A$ and $D$ equal
\[( 0,    -,     +,     +,     +,    -,     +,     +,    -,     +,     +,     +,    -)\quad\mbox{and}\quad
( -,     +,    -,    -,     +,     +,     +,     +,     +,     +,    -,    -,    +),\]
we get that $N^{-1} \in \bZ^{13 \times 13}$. Let $A_7, \ldots,A_{12}$, $S_7, \ldots, S_{12}$, and $B_7,\ldots, B_{12}$
be the corresponding $A$, $S=D|D|^{-1}$, $B_-=\sqrt{I-(1/5)A}$. We know from Corollary~\ref{Cor 1}
that $S_jB_j$ is a basis matrix for the $j$th lattice. It turns out that
$\det(5I-A_j)=2\,560\,000=2^{12} \cdot 5^4$ and hence again $\det B_j=2^6/5^{9/2} \approx 0.0458$ for $7 \le j \le 12$.

\medskip
Actually,
\begin{eqnarray*}
& & B_1=B_2, \quad B_3=B_4, \quad B_5=B_6,\\
& & S_7 B_7=-S_8 B_8, \quad S_9B_9=-S_{10}B_{10}, \quad  S_{11}B_{11}=-S_{12}B_{12},\\
& & B_1=U_1S_{11}B_{11}, \quad B_3=U_2 S_9B_9, \quad B_5=U_3S_7 B_7
\end{eqnarray*}
with orthogonal matrices $U_1, U_2, U_3$.
The relation ``$X \sim Y$ if and only if $XY^{-1}$ is a nonzero scalar multiple of an orthogonal matrix'' is an equivalence relation on every family of
invertible $k \times k$ matrices. The equivalence classes of this relation on
$\{B_1, \ldots, S_{12}B_{12}\}$ are
\[\{B_1=B_2, S_{11}B_{11}, S_{12}B_{12}\},\quad \{B_3=B_4, S_9B_9, S_{10}B_{10}\}, \quad
\{B_5=B_6, S_7 B_7, S_8 B_8\}.\]

The first rows of $(A_1,D_1), (A_3,D_3), (A_5,D_5)$ are
\begin{eqnarray*}
& & ( 0,    -,    -,    -,     +,    -,     +,     +,    -,     +,    -,    -,    -, \quad   -,    -,     +,     +,     +,    -,     +,     +,    -,     +,     +,   +,   -),\\
& & ( 0,     -,     +,     +,    -,    -,    -,    -,    -,    -,     +,     +,    -, \quad    +,    -,    -,    -,     +,    -,     +,     +,    -,     +,    -,   -, -),\\
& & ( 0,    +,     -,    -,    -,     +,    -,    -,     +,    -,    -,    -,     +,  \quad   +,    -,     +,     +,    -,    -,    -,    -,    -,    -,     +,   +,   -).
\end{eqnarray*}

\medskip
For $1 \le j \le 6$, the smallest eigenvalue of $B_j$ is about $0.3736$, whence
\[\n B_jx\n^2 > 0.37^2 \n x\n^2 > 0.13 \n x\n^2.\]
Thus, $\n Bx\n^2 >1$ for $\n x\n^2 \ge 7$. In the last $6$ cases, the smallest eigenvalue of $B_j$ is about
$0.4991$ and so we have
\[\n S_jB_jc\n^2=\n B_jx\n^2 > 0.49^2 \n x\n^2 > 0.24 \n x\n^2,\]
which is greater than $1$ for $\n x\n^2 \ge 4$. We took all $j \in \{1, \ldots, 12\}$ and $x \in \bZ^{13}$ with $\n x\n^2 \le 6$
and checked wether $\n B_jx\n^2 < 1.1$. For each $j$, we obtained
exactly $52$ vectors $x \in \bZ^{13}\sm \{0\}$ such that $\n B_jx\n^2 < 1.1$. The columns $B_jx$ are $\pm$
the $26$ columns $f_1, \ldots, f_{26}$ of $G$. Consequently, {\em in all cases the minimal distance of $\La(G)$ is $1$, $\La(G)$ has a basis
of minimal vectors, and $S(\La(G))=\{\pm f_1, \ldots, \pm f_{26}\}$.}
$\;\:\square$

\bigskip
{\bf Ten lattices from  (25,50) frames.} We finally take $k=25$, $n=50$, $\al=7$. We consider the $25 \times 25$ circulant matrices
$A$ and $D$ whose first rows are
\[(0,\eps_1,\eps_2,\eps_3,\eps_4,\eps_5,\eps_6, \eps_7,\eps_8,\eps_9,\eps_{10},\eps_{11},\eps_{12},
\eps_{12},\eps_{11},\eps_{10},\eps_9,\eps_8,\eps_7, \eps_6,\eps_5,\eps_4,\eps_3,\eps_2,\eps_1)\]
and
\begin{eqnarray*}
& & (\eps_{25},\eps_{13},\eps_{14},\eps_{15},\eps_{16},\eps_{17},\eps_{18}, \eps_{19},\eps_{20},\eps_{21},\eps_{22},\eps_{23},\eps_{24},\\
& & \qquad \:\eps_{24},\eps_{23},\eps_{22},\eps_{21},\eps_{20},\eps_{19}, \eps_{18},\eps_{17},\eps_{16},\eps_{15},\eps_{14},\eps_{13}),
\end{eqnarray*}
with $\eps_k \in \{-1,1\}=:\{-,+\}$, respectively.  These are $2^{25}=33\,554\,432$ matrices.
In exactly $20$ cases the matrix $C$ given by~(\ref{CZirk}) satisfies $C^2=49\,I$. One such case is
where the first rows of $A$ and $D$ are
\[(0,    -,    -,    -,     +,    -,     +,     +,    -,     +,     +,     +,
-,    -,     +,     +,     +,    -,     +,     +,    -,     +,    -,    -,   -)\]
and
\[(-,    -,     +,     +,     +,     +,     +,    -,     +,    -,     +,     +,
-,    -,     +,     +,    -,     +,    -,     +,     +,     +,     +,     +,   -),\]
respectively. We have
\[|\det D|=\det(7I+A)=\det(7I-A)=260\,119\,840^2=2^{22}\cdot 3^2\cdot 5^2\cdot 7^2\cdot 11^4,\]
$N \in \bZ^{25 \times 25}$ and $N^{-1} \in \bZ^{25 \times 25}$ in each of the $20$ cases.
Thus, by Corollary \ref{Cor 1}, we obtain $20$ lattices $\La(G_j)$ with $B_j=\sqrt{I+(1/7)A_j}$
as a basis matrix and
\[\det B_j=\frac{2^{11}\cdot 3 \cdot 5 \cdot 7\cdot 11^2}{7^{25/2}}
=\frac{2^{11}\cdot 3 \cdot 5 \cdot 11^2}{7^{23/2}} \approx 0.0007\,1052\]
for all $1 \le j \le 20$. In fact $B_j=B_{j+10}$ for $1 \le j \le 10$, and the equivalence classes of
the set $\{B_1, \ldots, B_{10}\}$
under the equivalence relation
``$B_i \sim B_j$ if and only if $B_iB_j^{-1}$ is a nonzero scalar multiple of an orthogonal matrix''
are the ten singletons $\{B_1\}, \ldots, \{B_{10}\}$.
The smallest eigenvalue of $B_j$ is about $0.1415$ for all $j$. $\;\:\square$

\section{Unit tight equiangular $\boldsymbol{(k,k+1)}$ frames}

Sometimes it is advantageous to represent a unit tight equiangular $(k,n)$ frame by coordinates
different from those in $\bR^k$. This is in particular the case for $(k+1,k)$ frames.

\medskip
Fix $k \ge 2$ and consider the set $\cF$ of the $k+1$ normalized columns of height $k+1$ formed by the permutations
of $-k, 1, \ldots, 1$ ($k$ ones),
\[f_1=\frac{1}{\sqrt{k^2+k}}\left(\begin{array}{r}
-k\\ 1\\ \vdots \\1\end{array}\right), \;
f_2=\frac{1}{\sqrt{k^2+k}}\left(\begin{array}{r}
1\\ -k\\ \vdots \\1\end{array}\right), \; \ldots, \;
f_{k+1}=\frac{1}{\sqrt{k^2+k}}\left(\begin{array}{r}
1\\ 1\\ \vdots \\-k\end{array}\right).\]
These $k+1$ vectors are in the orthogonal complement of $(1,\ldots,1)' \in \bR^{k+1}$ and may
therefore be thought of as vectors in $\bR^k$. Let
\[\La(\cF)={\rm span}_\bZ \{f_1, \ldots, f_{k+1}\} \subset \bR^k.\]
The following theorem in conjunction with Proposition~\ref{eutaxy} proves Theorem~\ref{Theo I2}(a).

\begin{thm} \label{Theo 5}
The vectors $f_1, \ldots, f_{k+1}$ form a unit tight equiangular $(k,k+1)$ frame and $\La(\cF)$ is a full rank lattice.
The matrix $B$ constituted by $f_1, \ldots, f_k$,
\begin{equation}
B=\frac{1}{\sqrt{k^2+k}}\left(\begin{array}{rrrr}
-k & 1 & \ldots & 1\\
1 & -k & \ldots & 1\\
\vdots & \vdots & & \vdots\\
1 & 1 & \ldots & -k\\
1 & 1 & \ldots &1\end{array}\right)_{(k+1) \times k},\label{4}
\end{equation}
is a basis matrix for $\La(\cF)$, we have
\[\det(B'B)=\frac{1}{{k+1}}\left(1+\frac{1}{k}\right)^{k},\]
the lattice $\La(\cF)$ has a basis of minimal vectors, and $S(\La(\cF))=\{\pm f_1, \ldots,\pm f_{k+1}\}$.
\end{thm}

{\em Proof.}
It is well known that $\cF$ is a tight unit equiangular $(k,k+1)$ frame.
We include the proof for the reader's convenience.
First, the columns of the matrix $B$ are easily seen to be linearly independent,
which shows that ${\rm span}_\bR \{f_1, \ldots, f_k\}=\bR^k$.
Secondly, it is clear that $\n f_j\n =1$ for all $j$.
Thirdly, we have $(f_i,f_j)=(-k-1)/(k^2+k)=-1/k$ for $i \neq j$.
And finally, if $x=(x_1, \ldots, x_{k+1})$ and $x_1+\cdots+x_{k+1}=0$, then
\[(f_j,x)=\frac{1}{\sqrt{k^2+k}}\left(-kx_j+\sum_{i \neq j}x_i\right)
= \frac{1}{\sqrt{k^2+k}}(-kx_j-x_j)\]
and hence
\[\sum_{j=1}^{k+1}(f_j,x)^2=\frac{1}{k^2+k}\sum_{j=1}^{k+1} (-(k+1)x_j)^2=\frac{k+1}{k}\n x\n^2,\]
that is, the frame is tight with $\ga=(k+1)/k$.

\medskip
Since $f_1+\cdots+f_k=-f_{k+1}$, we have $\La(\cF)={\rm span}_\bZ\{f_1, \ldots, f_k\}$. This shows
that $\La(\cF)$ is $\{BX: X \in \bZ^{k}\}$. Consequently, $\La(\cF)$ is a full rank lattice with the
matrix $B$ given by (\ref{4}) as a basis matrix. The product $B' B$ is
\begin{equation}
B' B=\frac{1}{k^2+k}\left(\begin{array}{rrrr}
a & b & \ldots & b\\
b & a & \ldots & b\\
\vdots & \vdots & \ddots & \vdots\\
b & b & \ldots & a \end{array}\right)_{k \times k}\label{5}
\end{equation}
with $a=k^2+k$ and $b=-k-1$. The determinant of a matrix of the form (\ref{5}) is known to
be $(a-b)^{k-1}(a+(k-1)b)$. Thus,
\[\det B' B=\frac{1}{(k^2+k)^k}(k^2+k+k+1)^{k-1}(k^2+k-(k-1)(k+1))=\frac{(k+1)^{k-1}}{k^k}.\]

\medskip
We are left with determining $S(\La(\cF))$. Straightforward computation shows that the inequality $\n Bx\n^2 \ge 1$ is equivalent to the inequality
\begin{equation}
(k+1)(x_1^2+\cdots+x_k^2) \ge k+(x_1+\cdots+x_k)^2,\label{ieq}
\end{equation}
and that equality holds in both inequalities only simultaneously. We first show (\ref{ieq}) for integers $(x_1, \ldots,x_k) \in \bZ^k \sm \{0\}$ by induction on $k$.
For $k=1$, inequality (\ref{ieq}) is trivial. Suppose it is true for $k-1$:
\[k(x_1^2+\cdots+x_{k-1}^2) \ge k-1+(x_1+\cdots +x_{k-1})^2.\]
If $x_1^2+\cdots+x_{k-1}^2 \ge 1$, we may add $x_1^2+\cdots+x_{k-1}^2$ on the left and $1$ on the right to get
\[(k+1)(x_1^2+\cdots+x_{k-1}^2) \ge k+(x_1+\cdots +x_{k-1})^2.\]
This proves (\ref{ieq}) in the case where one of the integers $x_1, \ldots, x_k$ is zero and one of them is nonzero.
We are so left with the case where $x_j \neq 0$ for all $j$. Then $x_1^2+\cdots+x_k^2 \ge k$ and hence
\begin{eqnarray}
& & k + (x_1+\cdots+x_k)^2 \le k+ (|x_1|+\cdots +|x_k|)^2 \label{i1} \\
& & \le k +k(x_1^2+\cdots+x_k^2) \nonumber \\
& & \le x_1^2+\cdots+x_k^2 + k(x_1^2+\cdots+x_k^2)=(k+1)(x_1^2+\cdots+x_k^2), \label{i3}
\end{eqnarray}
which completes the proof of (\ref{ieq}). At this point we have shown that $\{f_1, \ldots, f_k\}$
is a basis of minimal vectors.

\medskip
To identify all of $S(\La(\cF))$, we have to check when equality in (\ref{ieq}) holds.
Suppose first that $x_j \neq 0$ for all $j$. In that case we have (\ref{i1}) to (\ref{i3}).
Equality in  (\ref{i3}) holds if and only if $|x_j|=1$ for all $j$, and equality in~(\ref{i1})
is valid if and only if all the $x_j$ have the same sign. Thus, we get the two vectors
$x=(1,\ldots,1)'$ and $x=(-1,\ldots,-1)'$. The corresponding products $Bx$ are $-f_{k+1}$
and $f_{k+1}$. Suppose finally that one of the $x_j$ is zero, say $x_k=0$. From (\ref{ieq}) with
$k$ replaced by $k-1$ we know that
\[k(x_1^2+\cdots+x_{k-1}^2) \ge k-1+(x_1+\cdots +x_{k-1})^2.\]
If $x_1^2+\cdots+x_{k-1}^2 >1$, we may add this inequality to the previous one to obtain that
\[(k+1)(x_1^2+\cdots+x_{k-1}^2) \ge k+(x_1+\cdots +x_{k-1})^2.\]
Consequently, for $x_1^2+\cdots+x_{k-1}^2 >1$ equality in (\ref{ieq}) does not hold.
If $x_1^2+\cdots+x_{k-1}^2 =1$, then $x_j = \pm 1$ for some $j$ and $x_i=0$ for all $i \neq j$.
In that case equality in~(\ref{ieq}) holds and the vector $Bx$ is $\pm f_j$. In summary, we
have proved that the set $S(\La(\cF))$ of all minimal vectors is just $\{\pm f_1, \ldots, \pm f_{k+1}\}$.
$\;\:\square$

\section{The remaining frames in dimensions at most 9}

Recall that (\ref{flist}) lists the unit tight equiangular frames in dimensions $k \le 9$
different from the $(k,k+1)$ frames.
By Proposition~\ref{Prop I1}, the $(3,6)$, $(7,14)$, and $(9,18)$ frames do not yield lattices, and the
lattices resulting from the $(5,10)$ case were discussed in Section~4. We are left with
the $(6,16)$ and $(7,28)$ cases.

\bigskip
{\bf A lattice from a (6,16) frame.}
In \cite{Sus} we see the unit tight equiangular $(6,16)$ frame
\[G=\frac{1}{\sqrt{6}}\left(\begin{array}{cccccccccccccccc}
 + & + & + & + & + & + & + & + & + & + & + & + & + & + & + & + \\
   + & + & + & + & + & + & + & + & - & - & - & - & - & - & - & - \\
   + & + & + & + & - & - & - & - & + & + & + & + & - & - & - & - \\
   + & + & - & - & + & + & - & - & + & + & - & - & + & + & - & - \\
   + & - & + & - & + & - & + & - & + & - & + & - & + & - & + & - \\
   + & - & - & + & - & + & + & - & - & + & + & - & + & - & - & + \end{array}\right).\]
Here $GG' = (16/6)I$ and $G' G=I+(1/3)C$ with a $16 \times 16$ matrix $C$
whose diagonal entries are zero and the other entries of which are $\pm 1$.
 The six columns $f_1,f_2,f_3,f_4,f_5, f_9$ of the matrix $G$ are linearly
independent and each of the remaining $10$ columns is a linear combination with integer coefficients
of these six columns. Consequently, by Proposition~\ref{Prop 1} with $\be=1$, these six columns form a basis matrix,
\[B=\frac{1}{\sqrt{6}}\left(\begin{array}{cccccc}
   + & + & + & + & + &     +     \\
   + & + & + & + & + &     -     \\
   + & + & + & + & - &     +     \\
   + & + & - & - & + &     +     \\
   + & - & + & - & + &     +     \\
   + & - & - & + & - &     -     \end{array}\right).\]
We have $\det (B'B)=2^6/3^6$.

\medskip
With $B' B=U'E U$, we get $\n B x\n^2=(EUx,Ux) \ge 0.48 \n x \n^2$, and this is at least $6$
if $\n x\n^2 \ge 13$. So consider the $x \in \bZ^6\sm \{0\}$ with $\n x\n^2 \le 13$. Such $x$ contain only $0,\pm 1, \pm 2, \pm 3$,
and using Matlab we checked that $\n Bx\n^2 <6.1$ for exactly $32$ nonzero $x$ of these $7^6-1= 117\,648$
possible $x$. The $32$ columns $Bx$ are just $\pm$ the columns of $\sqrt{6}G$. Thus, {\em $\La(\cF)$ has a basis of minimal vectors
and $S(\La(G))=\{\pm f_1, \ldots, \pm f_{16}\}$.} $\;\:\square$

\bigskip
{\bf A perfect lattice from a (7,28) frame.}
It is well known that the $\tbinom{8}{2}=28$ vectors resulting
from the columns $(-3,-3, 1,1,1,1,1,1)'$ by permuting the entries form a tight equiangular $(7,28)$ frame.
To be precise, let $\cF$ be the set of the vectors
\[f_1=\frac{1}{\sqrt{24}}\left(\begin{array}{r} -3\\ -3\\ 1\\1\\1\\1\\1\\1\end{array}\right),
\quad \ldots,\quad f_{28}=\frac{1}{\sqrt{24}}\left(\begin{array}{r} 1\\1\\1\\1\\1\\1\\-3\\-3\end{array}\right).
\]
These are unit vectors in $\bR^8$. They are all orthogonal to the vector $(1,1,1,1,1,1,1,1)'$, and after identifying the orthogonal complement of this
vector with $\bR^7$, we may think of $f_1, \ldots, f_{28}$ as unit vectors in $\bR^7$. We consider the set
\[\La(\cF)={\rm span}_\bZ\{f_1,\ldots,f_{28}\} \subset \bR^7.\]
The columns of the $8 \times 7$ matrix
\[B=\frac{1}{\sqrt{24}}\left(\begin{array}{rrrrrrr}
-3 & -3 & -3 & -3 & -3 & -3 & 1\\
-3 & 1 & 1 & 1 & 1 & 1 & 1 \\
1 & -3 & 1 & 1 & 1 & 1 & -3 \\
1 & 1 & -3 & 1 & 1 & 1 & 1 \\
1 & 1 & 1 & -3 & 1 & 1 & 1\\
1 & 1 & 1 & 1 & -3 & 1 & -3\\
1 & 1 & 1 & 1 & 1 &-3 & 1\\
1 & 1 & 1 & 1 & 1 & 1 & 1\end{array}\right)\]
are formed by $7$ of the above $28$ vectors. We denote these $7$ vectors by $f_1, \ldots,f_7$.
For the reader's convenience, we show that $\{f_1, \ldots,f_{28}\}$ {\em is a tight unit equiangular $(7,28)$ frame.}
The rank of the matrix $B$ is $7$, and hence ${\rm span}_\bR\{f_1,\ldots,f_{28}\} = \bR^7$.
Clearly, $\n f_j\n =1$ for all $j$.
We have $|(f_i,f_j)|=8/24=1/3$ for $i \neq j$ (equiangularity).
Finally, let $x=(x_1, \ldots, x_8) \in \bR^7$. Then $x_1+\cdots+x_8=0$ and hence
\[\sum_j x_j^2+\sum_{j \neq k} x_jx_k =0,\]
which implies that
\[2 \sum_{j \neq k} x_jx_k = - 2 \n x\n^2.\]
We have
\begin{eqnarray*}
\sum_{\ell=1}^{28} (f_\ell,x)^2 & = & \frac{1}{24}\sum_{j < k}\left(-3x_j-3x_k+\sum_{m\neq j,k} x_m\right)^2
= \frac{1}{24}\sum_{j < k}(-3x_j-3x_k-x_j-x_k)^2\\
& = & \frac{2}{3}\sum_{j < k}(x_j+x_k)^2
 =  \frac{1}{3}\sum_{j \neq k}(x_j+x_k)^2
 =  \frac{1}{3}\sum_{j \neq k}(x_j^2+2x_jx_k+x_k^2)\\
& = & \frac{1}{3}(14\n x\n^2-2\n x\n^2)
=4\n x\n^2.
\end{eqnarray*}
This proves the tightness with $\ga =4$ (which, as is should be, is just $n/k=28/7$).

\medskip
Straightforward inspection shows that each of the vectors $f_8, \ldots, f_{28}$
is a linear combination with integer coefficients of the vectors $f_1,\ldots, f_7$.
Consequently, $\La(\cF)$ {\em is a full rank lattice in} $\bR^7$, $\{f_1, \ldots, f_7\}$ is a basis of $\La(\cF)$, and $B$ is a basis matrix.
We have
\[B' B=
\frac{1}{24}\left(\begin{array}{rrrrrrr}
 24  &   8   &  8  &   8  &   8  &   8  &  -8\\
     8  &  24 &    8 &    8 &   8  &   8 &    8\\
     8  &   8  &  24  &   8  &   8  &   8 &   -8\\
     8  &   8   &  8   & 24   &  8   &  8  &  -8\\
     8   &  8    & 8    & 8   & 24    & 8   &  8\\
     8    & 8     &8     & 8   &  8   & 24   & -8\\
    -8     & 8    & -8    & -8  &   8  &  -8  &  24\end{array}\right),\]
and straightforward computation gives
\[\det B' B= \frac{2^{27}}{24^7}=\frac{2^6}{3^7}.\]

\medskip
We now prove that the minimal norm of $\La(\cF)$ is $1$. Let
\[\ti{B}=\sqrt{24}B=\left(\begin{array}{rrrrrrr}
-3 & -3 & -3 & -3 & -3 & -3 & 1\\
-3 & 1 & 1 & 1 & 1 & 1 & 1 \\
1 & -3 & 1 & 1 & 1 & 1 & -3 \\
1 & 1 & -3 & 1 & 1 & 1 & 1 \\
1 & 1 & 1 & -3 & 1 & 1 & 1\\
1 & 1 & 1 & 1 & -3 & 1 & -3\\
1 & 1 & 1 & 1 & 1 &-3 & 1\\
1 & 1 & 1 & 1 & 1 & 1 & 1\end{array}\right).\]
Take $x \in \bZ^7$ and consider $y=\ti{B}x \in \bZ^8$. We are interested
in the $x$ for which $\n y\n^2 \le 24$. With $s:=x_1+\cdots+x_7$, we have
\begin{eqnarray*}
& & y_1=-3s+4x_7, \quad y_3=s-4x_2-4x_7, \quad y_6=s-4x_5-4x_7, \quad y_8=s,\\
& & y_2=s-4x_1, \quad y_4=s-4x_3, \quad y_5=s-4x_4, \quad y_7=s-4x_6.
\end{eqnarray*}
It suffices to search for all $x \in \bZ^7$ with $s \ge 0$ and $y_1^2+\cdots+y_8^2 \le 24$. This is impossible for $y_8=s >5$.
So we may assume that $0 \le s \le 4$.

\medskip
Suppose first that $s=4$. We then must have $y_1^2+\cdots+y_7^2 \le 9$. Since $y_1$ is an even number, it cannot be $\pm 3$.
Consequently, $-2 \le -3s+4x_7=-12+4x_7 \le 2$, which gives $x_7=3$. Analogously, as $y_3$ is even, we get
$-2 \le s-4x_2-4x_7 = -8-4x_2 \le 2$, which yields $x_2=-2$. In the same way we obtain $x_5=-2$. Finally, the even number
$s-4x_1=4-4x_1$ is at least $-2$, which implies that $x_1 \le 1$. Equally, $x_3, x_4,x_6 \le 1$. It follows that
\[s=x_1 +\cdots+x_7 \le 1+1+1+1-2-2+3 =3 < 4=s,\]
which is a contradiction.

\medskip
Thus, we may restrict our search to $0 \le s \le 3$ and $y_1^2+\cdots+y_7^2 \le 24$. The inequality $-4 \le -3s+4x_7 \le 4$
gives
\[-4 \le 3s-4 \le 4x_7 \le 3s+4 \le 13,\]
whence $-1 \le x_7 \le 3$. These are $5$ possibilities. From $-4 \le s-4x_j \le 4$ we obtain that $-1 \le x_j \le 1$
for $j=1,3,4,6$, which is $3^4$ possibilities, and the inequality $-4 \le s-4x_j-4x_7 \le 4$ delivers
\[-16 \le s-4+4x_7 \le 4y_j \le 4+s-4x_7 \le 11\]
and hence $-4 \le x_j \le 2$ for $j=2,5$, leaving us with $7^2$ possibilities. In summary, we have to check
$5\cdot 3^4\cdot 7^2=19\,845$ possibilities. Matlab does this with integer arithmetics within a second. The result
is that $0 \le s \le 3$ and $y_1^2+\cdots+y_7^2 \le 24$ happens in exactly $50$ cases. One of these cases
is $y=0$, and in the remaining $49$ cases $y$ is $\pm$ one of the $2\cdot 28 = 56$ vectors $\sqrt{24}f_j$.
(Recall that, by symmetry, we restricted ourselves to $s \ge 0$. For $-3 \le s \le 3$ and $y_1^2+\cdots+y_7^2 \le 24$
to happen we would obtain exactly $57$ cases: the case $y=0$ and the $56$ vectors $y$ given by $\pm\sqrt{24}f_j$.)
This proves that {\em the minimal distance of $\La(G)$ is $1$, that $S(\La(G))=\{\pm f_1, \ldots,\pm f_{28}\}$,
and that $\La(G)$ has a basis of minimal vectors.} From Proposition~\ref{eutaxy} we deduce that {\em the lattice $\La(G)$ is strongly eutactic.}

\medskip
We finally show that {\em this $(7,28)$ frame generates a perfect lattice.}
We have shown that the
$28$ lattice vectors $f_1, \ldots, f_{28}$
are minimal vectors.  These vectors are given by their coordinates in the ambient $\bR^8$. We use a special  $7 \times 8$
matrix $A$ to transform these vectors isometrically into $\bR^7$. The $j$th row
of $A$ is
\[\frac{1}{\sqrt{j^2+j}} (1, \ldots, 1, -j, 0, \ldots,0)\]
with $j$ ones and $7-j$ zeros. We have $A=EA_0$ with $E={\rm diag}(1/\sqrt{j^2+j})_{j=1}^7$
and with $(1, \ldots,1,-j,0, \ldots,0)$ being the $j$th row of $A_0$.
We then get the $28$ minimal vectors $Af_j=EA_0f_j$ ($j=1, \ldots, 28$) in $\bR^7$.
These give us $28$ symmetric $7 \times 7$ matrices $C_j=E(A_0f_j)(A_0f_j)'E$. The lattice $\La(\cF)$ is perfect
if the real span of these $28$ matrices is the space of all $7 \times 7$ symmetric matrices. Each symmetric  $28 \times 28$
matrix may be written as $ETE$ with a symmetric matrix $T$, and hence we are left with showing
that each symmetric $28 \times 28$ symmetric matrix $T$ is a real linear combination of the matrices
$(A_0f_j)(A_0f_j)'$. For $k=1, \ldots,7$,
let \[([C_j]_{k,k}, [C_j]_{k+1,k}, \ldots, [C_j]_{7,k})'\] be the column formed by the entries of
the $k$th column of $C_j$ that are on or below the main diagonal. Stack these columns to a column $D_j$ of height
$7+6+\cdots+1=28$. The lattice is perfect if and only if the real span of $D_1, \ldots, D_{28}$ is all of
$\bR^{28}$, which happens if and only if the $28 \times 28$ matrix $D$ constituted by the $28$ columns
$D_1, \ldots, D_{28}$ is invertible. 
Tables~\ref{table2} and~\ref{table3} show the matrix $D$.

\begin{table}[h]
{\scriptsize
\caption{The first 14 columns of the matrix $D$.}
\label{table2}
\[\left(\begin{array}{rrrrrrrrrrrrrr}
 0 &16 &16 &16 &16 &16 &16 &16 &16 &16 &16 &16 &16 &0\\
     0 &-4 &12 &12 &12 &12 &12 &4 &-12 &-12 &-12 &-12 &-12 &0\\
     0 &24 &-8 &8 &8 &8 &8 &-24 &8 &-8 &-8 &-8 &-8 &0\\
     0 &20 &20 &-12 &4 &4 &4 &-20 &-20 &12 &-4 &-4 &-4 &0\\
     0 &16 &16 &16 &-16 &0 &0 &-16 &-16 &-16 &16 &0 &0 &0\\
     0 &12 &12 &12 &12 &-20 &-4 &-12 &-12 &-12 &-12 &20 &4 &0\\
     0 &8 &8 &8 &8 &8 &-24 &-8 &-8 &-8 &-8 &-8 &24 &0\\
    49 &1 &9 &9 &9 &9 &9 &1 &9 &9 &9 &9 &9 &25\\
    42 &-6 &-6 &6 &6 &6 &6 &-6 &-6 &6 &6 &6 &6 &10\\
    35 &-5 &15 &-9 &3 &3 &3 &-5 &15 &-9 &3 &3 &3 &-25\\
    28 &-4 &12 &12 &-12 &0 &0 &-4 &12 &12 &-12 &0 &0 &-20\\
    21 &-3 &9 &9 &9 &-15 &-3 &-3 &9 &9 &9 &-15 &-3 &-15\\
    14 &-2 &6 &6 &6 &6 &-18 &-2 &6 &6 &6 &6 &-18 &-10\\
    36 &36 &4 &4 &4 &4 &4 &36 &4 &4 &4 &4 &4 &4\\
    30 &30 &-10 &-6 &2 &2 &2 &30 &-10 &-6 &2 &2 &2 &-10\\
    24 &24 &-8 &8 &-8 &0 &0 &24 &-8 &8 &-8 &0 &0 &-8\\
    18 &18 &-6 &6 &6 &-10 &-2 &18 &-6 &6 &6 &-10 &-2 &-6\\
    12 &12 &-4 &4 &4 &4 &-12 &12 &-4 &4 &4 &4 &-12 &-4\\
    25 &25 &25 &9 &1 &1 &1 &25 &25 &9 &1 &1 &1 &25\\
    20 &20 &20 &-12 &-4 &0 &0 &20 &20 &-12 &-4 &0 &0 &20\\
    15 &15 &15 &-9 &3 &-5 &-1 &15 &15 &-9 &3 &-5 &-1 &15\\
    10 &10 &10 &-6 &2 &2 &-6 &10 &10 &-6 &2 &2 &-6 &10\\
    16 &16 &16 &16 &16 &0 &0 &16 &16 &16 &16 &0 &0 &16\\
    12 &12 &12 &12 &-12 &0 &0 &12 &12 &12 &-12 &0 &0 &12\\
     8 &8 &8 &8 &-8 &0 &0 &8 &8 &8 &-8 &0 &0 &8\\
     9 &9 &9 &9 &9 &25 &1 &9 &9 &9 &9 &25 &1 &9\\
     6 &6 &6 &6 &6 &-10 &6 &6 &6 &6 &6 &-10 &6 &6\\
     4 &4 &4 &4 &4 &4 &36 &4 &4 &4 &4 &4 &36 &4
\end{array}\right)\]}
\end{table}

\begin{table}
{\scriptsize
\caption{The last 14 columns of the matrix $D$.}
\label{table3}
\[\left(\begin{array}{rrrrrrrrrrrrrr}
    0 &0 &0 &0 &0 &0 &0 &0 &0 &0 &0 &0 &0 &0\\
     0 &0 &0 &0 &0 &0 &0 &0 &0 &0 &0 &0 &0 &0\\
     0 &0 &0 &0 &0 &0 &0 &0 &0 &0 &0 &0 &0 &0\\
     0 &0 &0 &0 &0 &0 &0 &0 &0 &0 &0 &0 &0 &0\\
     0 &0 &0 &0 &0 &0 &0 &0 &0 &0 &0 &0 &0 &0\\
     0 &0 &0 &0 &0 &0 &0 &0 &0 &0 &0 &0 &0 &0\\
     0 &0 &0 &0 &0 &0 &0 &0 &0 &0 &0 &0 &0 &0\\
    25 &25 &25 &25 &1 &1 &1 &1 &1 &1 &1 &1 &1 &1\\
   -10 &-10 &-10 &-10 &6 &6 &6 &6 &2 &2 &2 &2 &2 &2\\
    15 &-5 &-5 &-5 &3 &-1 &-1 &-1 &7 &7 &7 &3 &3 &3\\
   -20 &20 &0 &0 &-4 &4 &0 &0 &4 &0 &0 &8 &8 &4\\
   -15 &-15 &25 &5 &-3 &-3 &5 &1 &-3 &5 &1 &5 &1 &9\\
   -10 &-10 &-10 &30 &-2 &-2 &-2 &6 &-2 &-2 &6 &-2 &6 &6\\
     4 &4 &4 &4 &36 &36 &36 &36 &4 &4 &4 &4 &4 &4\\
    -6 &2 &2 &2 &18 &-6 &-6 &-6 &14 &14 &14 &6 &6 &6\\
     8 &-8 &0 &0 &-24 &24 &0 &0 &8 &0 &0 &16 &16 &8\\
     6 &6 &-10 &-2 &-18 &-18 &30 &6 &-6 &10 &2 &10 &2 &18\\
     4 &4 &4 &-12 &-12 &-12 &-12 &36 &-4 &-4 &12 &-4 &12 &12\\
     9 &1 &1 &1 &9 &1 &1 &1 &49 &49 &49 &9 &9 &9\\
   -12 &-4 &0 &0 &-12 &-4 &0 &0 &28 &0 &0 &24 &24 &12\\
    -9 &3 &-5 &-1 &-9 &3 &-5 &-1 &-21 &35 &7 &15 &3 &27\\
    -6 &2 &2 &-6 &-6 &2 &2 &-6 &-14 &-14 &42 &-6 &18 &18\\
    16 &16 &0 &0 &16 &16 &0 &0 &16 &0 &0 &64 &64 &16\\
    12 &-12 &0 &0 &12 &-12 &0 &0 &-12 &0 &0 &40 &8 &36\\
     8 &-8 &0 &0 &8 &-8 &0 &0 &-8 &0 &0 &-16 &48 &24\\
     9 &9 &25 &1 &9 &9 &25 &1 &9 &25 &1 &25 &1 &81\\
     6 &6 &-10 &6 &6 &6 &-10 &6 &6 &-10 &6 &-10 &6 &54\\
     4 &4 &4 &36 &4 &4 &4 &36 &4 &4 &36 &4 &36 &36
\end{array}\right)\]}
\end{table}

The matrix $D$ can be constructed with integer arithmetics. 
The determinant
$\det D$ may be computed by the Gaussian algorithm and thus with integer arithmetics, too.
In the intermediate steps, one may factor out powers of $2$. For example, in the original
matrix $D$ we may draw out $16$ from the first line, $4$ from the second, $8$ from the third,
and so on. It results that
\[\det D=16^3\cdot 8^4\cdot 4^9\cdot 2^6 \cdot \det \ti{D}= 2^{48} \det \ti{D},\]
and we may start the Gaussian algorithm with $\det\ti{ D}$.  The final result is \[\det D=3\cdot 2^{159}.\]
As this is nonzero, we conclude that $D$ is invertible and thus that $\La(\cF)$ is perfect. At this point the proof
of Theorem~\ref{Theo I2}(c) is complete.

\medskip
The perfection of this lattice was also established by  Bacher in~\cite{bacher} (see Section 7, especially 7.1). However, Bacher's approach is different from ours:
he obtains the lattice in question as the kernel of a certain linear map, establishes its perfection, and then remarks that its set of minimal vectors comprises an equiangular system.
We, on the other hand, construct the lattice from the equiangular frame and show its perfection directly from this construction. Hence our argument here complements Bacher's, going in the opposite direction.

\medskip
Since $\La(\cF)$ is perfect and strongly eutactic, the packing density of this lattice is a local maximum. As we know the
minimal distance and the determinant of this lattice, the packing density can be easily computed using~(\ref{den}).
It turns out to be $21.57\,\%$. This is better than the packing density of the root lattice $A_7$,
which is $14.76\,\%$. In~\cite{SIAM}, we studied lattices in $\bR^k$ that are generated by Abelian groups of the order $k+1$.
There the packing density of the lattices generated by Abelian groups of order $8$
was shown to $20.88\,\%$. Thus, $\La(\cF)$ is also better than this. We nevertheless do not reach the best
packing density for a $7$-dimensional lattice, which is $29.53\,\%$ and is achieved for the well known lattice $E_7$.
$\;\:\square$

\medskip

{\bf Acknowledgement.} Fukshansky acknowledges support of the NSA grant H98230-1510051. Garcia acknowledges support of the NSF grant DMS-1265973. Needell acknowledges support of the Alfred P. Sloan Fellowship and NSF Career grant number 1348721.

\bigskip
\bigskip
A. B\"ottcher, Fakult\"at f\"ur Mathematik, TU Chemnitz, 09107 Chemnitz, Germany

{\tt aboettch@mathematik.tu-chemnitz.de}

\medskip
L. Fukshansky, Department of Mathematics,  Claremont McKenna College,

850 Columbia Ave,
Claremont, CA 91711, USA

{\tt lenny@cmc.edu}

\medskip
S. R. Garcia, Department of Mathematics, Pomona College,

610 N. College Ave, Claremont, CA 91711, USA

{\tt stephan.garcia@pomona.edu}

\medskip
H. Maharaj, Department of Mathematics, Pomona College,

610 N. College Ave, Claremont, CA 91711, USA

{\tt hirenmaharaj@gmail.com}

\medskip
D. Needell, Department of Mathematics,  Claremont McKenna College,

850 Columbia Ave,
Claremont, CA 91711, USA

{\tt dneedell@cmc.edu}


\begin{thebibliography}{20}

\bibitem{axell2012spectrum}

E. Axell, G. Leus, E. G. Larsson, and H. V. Poor,
{\em Spectrum sensing for cognitive radio: State-of-the-art and recent
  advances.}
IEEE Signal Proc. Mag. 29(3) (2012),~101--116.

\bibitem{bacher} R. Bacher,
{\em Constructions of some perfect integral lattices with minimum 4.}
J. Th\'eor. Nombres Bordeaux 27 (2015), no. 3, 655--687.

\bibitem{baraniuk2008simple}
R. Baraniuk, M. Davenport, R. DeVore, and M. Wakin,
{\em A simple proof of the restricted isometry property for random
  matrices.}
Constr. Approx. 28 (2008), 253--263.

\bibitem{SIAM} A. B\"ottcher, L. Fukshansky, S. R. Garcia, and H. Maharaj,
{\em On lattices generated by finite Abelian groups.}
SIAM J. Discrete Math. 29 (2015), 382--404.

\bibitem{candes2005error}
E. Candes, M. Rudelson, T. Tao, and R. Vershynin,
{\em Error correction via linear programming.}
In: IEEE Symp. Found. Comput. Sci.
  (FOCS'05),  IEEE, 2005, 668--681.
	
\bibitem{fickus2012steiner}
M. Fickus, D. G. Mixon, and J. C. Tremain,
{\em Steiner equiangular tight frames.}
Linear Algebra Appl. 436 (2012), 1014--1027.


\bibitem{flinth2016promp}
A. Flinth and G. Kutyniok,
{\em Promp: A sparse recovery approach to lattice-valued signals.}
Preprint 2016.

\bibitem{Goe} J. M. Goethals and J. J. Seidel,
{\em Orthogonal matrices with zero diagonal.} Canad. J. Math. 19 (1967), 1001--1010.

\bibitem{HP} R. B. Holmes and V. Paulsen,
{\em Optimal frames for erasures.} Linear Algebra Appl. 377 (2004), 31--51.


\bibitem{martinet}
J.~Martinet,
\newblock {\em Perfect Lattices in Euclidean Spaces}.
\newblock Springer-Verlag, 2003.

\bibitem{mixon2011equiangular}
D. G. Mixon, C. Quinn, N. Kiyavash, and M. Fickus,
{\em Equiangular tight frame fingerprinting codes.}
In: IEEE Int. Conf. Acoustics, Speech and
  Sig. Proc. (ICASSP), IEEE, 2011, 1856--1859.


\bibitem{Pa} R. E. A. C. Paley, {\em On orthogonal matrices.} J. Math. Phys. 12 (1933), 311--320.

\bibitem{rossi2014spatial}
M. Rossi, A. M. Haimovich, and Y. C. Eldar,
{\em Spatial compressive sensing for mimo radar.}
IEEE Trans. Signal Process., 62(2) (2014), 419--430.


\bibitem{rusu2015optimized}
C. Rusu and N. Gonz{\'a}lez-Prelcic,
{\em Designing incoherent frames through convex techniques for optimized compressed sensing.}
IEEE Trans. Signal Process. 64(9) (2016), 2334--2344.


\bibitem{achill1} A. Sch\"urmann, {\em Perfect, strongly eutactic lattices are periodic extreme.} Adv. Math. 225 (2010), no. 5, 2546--2564.

\bibitem{strohmer2003optimal}
T. Strohmer and S. Beaver,
{\em Optimal {OFDM} design for time-frequency dispersive channels.}
IEEE Trans. Commun., 51(7) (2003), 1111--1122.


\bibitem{Stroh} T. Strohmer and R. W. Heath Jr., {\em Grassmannian frames with applications to coding and communication.}
Appl. Comput. Harmon. Anal. 14 (2003), 257--275.


\bibitem{Sus} M. A. Sustik, J. A. Tropp, I. S. Dhillon, and R. W. Heath Jr.,
{\em On the existence of equiangular tight frames.} Linear Algebra Appl. 426 (2007), 619--635.

\bibitem{tropp2004greed}
J. A. Tropp,
{\em Greed is good: Algorithmic results for sparse approximation.}
IEEE Trans. Inform. Theory, 50(10) (2004), 2231--2242.


\bibitem{tsiligianni2012use}
E. Tsiligianni, L. P. Kondi, and A. K. Katsaggelos,
{\em Use of tight frames for optimized compressed sensing.}
In: Proc. Signal Process. Conf. (EUSIPCO), IEEE, 2012, 1439--1443,

\end{thebibliography}
\end{document}